\documentclass[letterpaper,oneside,10pt]{article}

\usepackage[USenglish]{babel} 
\usepackage[T1]{fontenc}
\usepackage[ansinew]{inputenc}

\usepackage{lmodern} 

\usepackage{graphicx} 

\usepackage{amsmath}
\usepackage{amsthm}
\usepackage{amsfonts}
\usepackage[hyphens]{url}

\let\Right\right 
\let\Left\left 
\makeatletter 
\def\right#1{\Right#1\@ifnextchar){\!\right}{}} 
\def\left#1{\Left#1\@ifnextchar({\!\left}{}} 
\makeatother

\begin{document}

\pagestyle{empty} 


\title{Extended calculations of a special Harmonic number}
\author{John Blythe Dobson (j.dobson@uwinnipeg.ca)}

\maketitle


\pagestyle{plain} 


\begin{abstract}
\noindent
The search for values of $p$ for which the Harmonic numbers $H_{\lfloor p/6 \rfloor}$ vanish mod $p$, carried to $p < 600,000$ by Schwindt \cite{Schwindt} in 1983, is extended here to $p < 149,250,000,000,000$, and two new solutions are reported. (These results can now be found in the Online Encyclopedia of Integer Sequences, entry no.\ A238201.)

\noindent
\textit{Keywords}: Harmonic numbers, Fermat quotient, Fermat's Last Theorem
\end{abstract}

\noindent
In 1983 Schwindt \cite{Schwindt} sought solutions of the Harmonic number congruence

\begin{equation} \label{eq:Harmonic}
H_{\lfloor p/6 \rfloor} := \sum_{j=1}^{\lfloor p/6 \rfloor} \frac{1}{j} \equiv 0 \pmod{p},
\end{equation}

\noindent
with $\lfloor \cdot \rfloor$ denoting the greatest-integer function, and found the single case $p=61$ with $p < 600,000$. The three zeros found by us are:

\begin{displaymath}
61, \medspace 1\,680\,023, \medspace 7\,308\,036\,881.
\end{displaymath}

\noindent
We know of no intervening computations of these numbers \textit{per se}, though as we afterwards learned and as will be explained below, our solutions had already appeared in another guise. The historical motivation for the study of this sum is the famous proof given by Emma Lehmer in 1938 (\cite{Lehmer1938}, p.\ 358) that (\ref{eq:Harmonic}) is a necessary condition on the exponent $p$ for the failure of the first case of Fermat's Last Theorem (FLT). This result retains its interest despite the full proof of FLT by Andrew Wiles. A survey of related results appears in \cite{Ribenboim1979}, and subsequent progress in this direction has been made in \cite{Sun+Sun1992} and \cite{DilcherSkula}. It should however be noted that a statement in \cite{DilcherSkula}, pp.\ 389--390, implying that $H_{\lfloor p/N \rfloor} \equiv 0$ mod $p$ has a solution $p < 2,000$ for every $N$ between 2 and 46 other than 5, is incorrect (at least if we require $p > N$ to ensure that the sum is not vacuous). In fact, we have determined that there are no solutions with $p < 319,900,000$ for $N = 5, 12, 17, 18, 20, 29, 31, 43$, and for $N = 5, 12$ there are no solutions with $p < 9,520,000,000,000$ (see \cite{DobsonCalculationsSpecialHarmonicNumbers}).

Schwindt reports obtaining the right-hand side of (\ref{eq:Harmonic}) by direct evaluation, and states that this operation constituted nearly half the work for a project described as having been ``run at night over half a year.'' As he cites Lehmer's work, it is not clear why he chose to perform the calculations in such a processing-intensive manner, rather than to exploit her result that

\begin{equation} \label{eq:H6}
H_{\lfloor p/6 \rfloor} \equiv -2 \cdot q_p(2)  - \frac{3}{2} \cdot q_p(3) \quad \pmod{p},
\end{equation}

\noindent
where the $q_p(b) := (b^{p - 1} - 1)/p$ are Fermat quotients. On Schwindt's own showing, his method has an algorithmetic complexity of order $p$ log $p$, because the calculation of each modular inverse has order log $p$, while the range over which they must be summed has order $p$. In contrast, the complexity of each of the Fermat quotients in (\ref{eq:H6}) has only order log $p$, the same as that of a single modular inverse, and it is only necessary to calculate two of them. (In fact, as will be shown below, it is possible to treat the two together; and the ensuing analysis assumes this saving.) In other words, over a range of $p$ running from some small number $\epsilon$ to $n$, the processing of Schwindt's calculations would have order

\begin{displaymath} \label{eq:SchwindtCumulativeTime}
\int_\epsilon^n x\,\mathrm{log}\,x\,\mathrm{d}x \sim \frac{n^2}{2}\left(\mathrm{log}\,n - \frac{1}{2}\right),
\end{displaymath}

\noindent
while that for the right-hand side of (\ref{eq:H6}) would have only order

\begin{displaymath} \label{eq:LehmerCumulativeTime}
\int_\epsilon^n \mathrm{log}\,x\,\mathrm{d}x \sim n(\mathrm{log}\,n - 1),
\end{displaymath}

\noindent
and the results of our tests are in close agreement with these predictions. These models fail to take into account the diminishing frequency of the primes as $n$ increases and thus are not fully realistic, but asymptotically the effect of any further refinement would be dwarfed by the fundamental differences between the two methods.

In Table 1 below, we give a comparison of runtimes for the Fermat-quotient method versus that of Schwindt, for various upper limits $n$: 600,000 (the limit of Schwindt's calculations), 1,680,023 (the first large zero), and 7,308,036,881 (the second large zero). The last limit is attainable only in the case of the Fermat-quotient test, not Schwindt's. These calculations were performed in Mathematica on a typical desktop computer, with a 3.2 GHz processor and 4 GB of memory, with the network connection disabled to prevent interruptions. Each run was begun at the number 7 (the least prime for which $H_{\lfloor p/6 \rfloor}$ is meaningfully defined) and extended precisely to the stated value of the limit $n$. Though runtimes inevitably vary over the course of multiple trials, these were performed under consistent conditions and may be assumed to give a good relative sense of the times involved, so we have reported them precisely as obtained, rounded to the nearest second.

\begin{table} [hb]
\begin{center}
\caption{Comparison of timings for calculations of zeros of $H_{\lfloor p/6 \rfloor}$ by various methods}
\label{Table_1}
\begin{tabular}{ r | r | r }
upper limit $n$  & Fermat-quotient test & Schwindt's method \\
\hline
600,000        &           2 seconds &   1,633 seconds \\
1,680,023      &           5 seconds &  12,422 seconds \\
7,308,036,881  &      39,864 seconds &    (impossible) \\
\end{tabular}
\end{center}
\end{table}

These results suggest that Schwindt's calculations up to the limit $n = 600\,000$ could not have been profitably pursued much further using the equipment then available. Indeed, even with today's processing capabilites, such calculations could only be extended about a hundredfold. It would be inconceivable to find the third zero using Schwindt's test, as the estimated runtime is about 10,000 years; but this objection is moot because even an attempt to spot-check this solution by Schwindt's method produces a memory allocation failure error in Mathematica. In contrast, the Fermat-quotient test obtains the result in about 39,864 seconds (just over 11 hours).

We shall now describe how the calculation based on Emma Lehmer's congruence (\ref{eq:H6}) was optimized. Because we are only interested in locating its zeros and do not require the actual values, we are free to disregard sign and scale, and may multiply throughout by $-2$ to clear the negative signs and the fraction, establishing that for $p > 5$ the zeros coincide with those of

\begin{displaymath}
4 \cdot q_p(2)  + 3 \cdot q_p(3).
\end{displaymath}

\noindent
Applying in reverse the logarithmetic and factorization rules for the Fermat quotient given by Eisenstein \cite{Eisenstein}, this expression may be realized more compactly as follows:

\begin{displaymath}
q_p(2^4) + q_p(3^3) \equiv q_p(2^4 \cdot 3^3) \equiv q_p(432) \quad \pmod{p} \quad (p > 5).
\end{displaymath}

\noindent
By simplification and consolidation of the Fermat quotient calculation, we are able to improve the runtime considerably. Moreover, this device alerted us to the fact that pertinent prior literature might include studies of the divisibility of the Fermat quotient by $p$ for composite bases, and such a study was indeed located in the form of Le\v{z}\'{a}k \cite{Lezak}, who performed similar calculations to the limit $p < 35,000,000,000$ and whose results agree perfectly with out own. (We regret having inadvertently deprived Le\v{z}\'{a}k of due credit in previous versions of these notes.) The more extended calculations begun in 2014 by Richard Fischer \cite{Fischer} have since reached the limit $p < 163,000,000,000,000$ as of 24 December 2022, without finding further solutions in this case. This limit is six orders of magnitude higher than any calculation for $H_{\lfloor p/6 \rfloor}$ that could be achieved by direct calculations of sums of reciprocals in under a decade of processing time. We suspended our own calculations on learning of Fischer's, and have independently confirmed them onlly to $p < 149,250,000,000,000$.

\section*{Acknowledgement}

I am grateful to T.\,D.\ Noe for creating OEIS sequence A238201 based on an earlier version of this report.

\end{document}